\documentclass{amsart}
\usepackage{amsmath,amssymb,amsxtra,amsthm,amscd}
\usepackage[mathscr]{eucal}
\usepackage{bbm}
\usepackage{graphicx}
\usepackage[all,cmtip]{xy}
\usepackage{rotating}
\usepackage{lmodern}
\usepackage{rotating}
\usepackage{microtype}

\newif\ifShowLabels
\ShowLabelstrue
\newcommand{\TeXref}[1]{
\marginpar{\scriptsize \texttt{#1}}}
\ShowLabelsfalse

\let\oldtocsection=\tocsection
\let\oldtocsubsection=\tocsubsection
\renewcommand{\tocsection}[2]{\hspace{0em}\oldtocsection{#1}{#2}}
\renewcommand{\tocsubsection}[2]{\hspace{2em}\oldtocsubsection{#1}{#2}}

\DeclareMathOperator{\B}{\mathbf{B}}

         \newcommand{\Btimes}{\B_{\times}}

\DeclareMathOperator{\CAT}{CAT}

\DeclareMathOperator{\EGamma}{\boldsymbol{E}\boldsymbol{\Gamma}}
\DeclareMathOperator{\End}{End}

\DeclareMathOperator{\Free}{\mathbf{Free}}

\DeclareMathOperator{\Fun}{Fun}

         \newcommand{\Gnc}{G
         }

                  \newcommand{\Itimes}{{I}_{\times}}

\DeclareMathOperator{\Hom}{Hom}

\DeclareMathOperator{\id}{id}

\DeclareMathOperator{\LI}{\mathbf{LI}}

\DeclareMathOperator{\Mod}{\mathbf{Mod}}

         \newcommand{\subdot}{\boldsymbol{\cdot}}

\DeclareMathOperator*{\one}{1}
\newcommand{\onehatplace}[1]
{ \one^{\substack{#1 \\ \frown}} }

\DeclareMathOperator*{\bones}{\times}
\newcommand{\undertimes}[1]
{ \bones_{#1} }

\DeclareMathOperator*{\bowl}{\cup}
\newcommand{\undercup}[1]
{ \bowl_{#1} }

\DeclareMathOperator*{\arch}{\cap}
\newcommand{\undercap}[1]
{ \arch_{#1} }

\newcommand{\pull}
{\!\!\! -\!\!\! -\!\!\! -\!\!\!}

\DeclareMathOperator*{\holimprep}{holim}                       
\newcommand{\holim}[1]%
{\displaystyle\holimprep_{\substack{\leftarrow \pull - \\ #1}} \, }

\DeclareMathOperator*{\hocolimprep}{hocolim}                   
\newcommand{\hocolim}[1]%
{\displaystyle\hocolimprep_{\substack{- \pull \rightarrow \\ #1}} \, }

\DeclareMathOperator*{\plainlim}{lim}                           
\newcommand{\contralim}[1]%
{\displaystyle\plainlim_{\substack{\leftarrow \pull - \\ #1}} \, }

\DeclareMathOperator*{\plaincolim}{colim}                       
\newcommand{\colim}[1]%
{\displaystyle\plaincolim_{\substack{- \pull \rightarrow \\ #1}} \, }

\DeclareMathOperator*{\laxlimplain}{laxlim}                     
\newcommand{\laxlim}[1]%
{\displaystyle\laxlimplain_{\substack{\leftarrow \pull - \\ #1}} \, }

\providecommand{\bysame}{\makebox[3em]{\hrulefill}\thinspace}

\theoremstyle{plain}
\newtheorem{Thm}{Theorem}[section]

\newtheorem*{MainThm}{Main Theorem}

\newtheorem{Cor}[Thm]{Corollary}

\newtheorem{Lem}[Thm]{Lemma}
\newtheorem{Prop}[Thm]{Proposition}

\theoremstyle{definition}
\newtheorem{Def}[Thm]{Definition}

\newtheorem{Ex}[Thm]{Example}

\newtheorem{Rem}[Thm]{Remark}

\theoremstyle{remark}
\newtheorem{Not}[Thm]{Notation}

\newtheoremstyle{freestylethm}{6pt}{6pt}{\itshape}{}%
                {\bfseries}{}{.5em}{\thmnote{#3}}
\theoremstyle{freestylethm}

\newcommand{\SecRef}[2]{\section{#1}\label{S:#2}%
\ifShowLabels \TeXref{{S:#2}} \fi}

\newcommand{\refT}[1]{\textup{\ref{T:#1}}}
\newcommand{\refL}[1]{\textup{\ref{L:#1}}}
\newcommand{\refD}[1]{\textup{\ref{D:#1}}}
\newcommand{\refC}[1]{\textup{\ref{C:#1}}}

\newenvironment{ThmRef}[1]%
{ \begin{Thm} \label{T:#1}
\ifShowLabels \TeXref{T:#1} \fi }%
{ \end{Thm} }
\newenvironment{DefRef}[1]%
{ \begin{Def} \label{D:#1}
\ifShowLabels \TeXref{D:#1} \fi }%
{ \end{Def} }
\newenvironment{LemRef}[1]%
{ \begin{Lem} \label{L:#1}
\ifShowLabels \TeXref{L:#1} \fi }%
{ \end{Lem} }
\newenvironment{CorRef}[1]%
{ \begin{Cor} \label{C:#1}
\ifShowLabels \TeXref{C:#1} \fi }%
{ \end{Cor} }
\newenvironment{RemRef}[1]%
{ \begin{Rem} \label{R:#1}
\ifShowLabels \TeXref{R:#1} \fi }%
{ \end{Rem} }
\newenvironment{PropRef}[1]%
{ \begin{Prop} \label{P:#1}
\ifShowLabels \TeXref{P:#1} \fi }%
{ \end{Prop} }
\newenvironment{ExRef}[1]%
{ \begin{Ex} \label{E:#1}
\ifShowLabels \TeXref{E:#1} \fi  }%
{ \end{Ex} }
\newenvironment{NotRef}[1]%
{ \begin{Not} \label{N:#1}
\ifShowLabels \TeXref{N:#1} \fi }%
{ \end{Not} }

\newenvironment{ThmRefName}[2]%
{ \begin{Thm} [#2]\label{T:#1}
\ifShowLabels \TeXref{T:#1} \fi }%
{ \end{Thm} }
{ \begin{Def} [#2]\label{D:#1}
\ifShowLabels \TeXref{D:#1} \fi }%
{ \end{Def} }
{ \begin{Lem} [#2]\label{L:#1}
\ifShowLabels \TeXref{L:#1} \fi }%
{ \end{Lem} }
{ \begin{Cor} [#2]\label{C:#1}
\ifShowLabels \TeXref{C:#1} \fi }
{ \end{Cor} }
{ \begin{Rem} [#2]\label{R:#1}
\ifShowLabels \TeXref{R:#1} \fi }%
{ \end{Rem} }
{ \begin{Prop} [#2]\label{P:#1}
\ifShowLabels \TeXref{P:#1} \fi }%
{ \end{Prop} }
{ \begin{Ex} [#2]\label{E:#1}
\ifShowLabels \TeXref{E:#1} \fi }%
{ \end{Ex} }

\begin{document}
\title{On Modules over Infinite Group Rings}
\author[Gunnar Carlsson]{Gunnar Carlsson}
\address{Department of Mathematics\\
Stanford University\\ Stanford\\ CA 94305}
\email{gunnar@math.stanford.edu}
\author[Boris Goldfarb]{Boris Goldfarb}
\address{Department of Mathematics and Statistics\\
SUNY\\ Albany\\ NY 12222} \email{goldfarb@math.albany.edu}
\thanks{The authors acknowledge support from the National
Science Foundation.}
\date{\today}

\begin{abstract}
Let $R$ be a commutative ring and $\Gamma$ be an infinite discrete group.
The algebraic $K$-theory of the group ring $R[\Gamma]$ is an important object of computation in geometric topology and number theory.  
When the group ring is Noetherian there is a companion $G$-theory of $R[\Gamma]$ which is often easier to compute.
However, it is exceptionally rare that the group ring is Noetherian for an infinite group.

In this paper, we define a version of $G$-theory for any finitely generated discrete group.
This construction is based on the coarse geometry of the group.
Therefore it has some expected properties such as independence from the choice of a word metric.
We prove that, whenever $R$ is a regular Noetherian ring of finite global homological dimension and $\Gamma$ has finite asymptotic dimension and a finite model for the classifying space $B\Gamma$, the natural Cartan map from the $K$-theory of $R[\Gamma]$ to $G$-theory is an equivalence.
On the other hand, our $G$-theory is indeed better suited for computation as we show in a separate paper.

Some results and constructions in this paper might be of independent interest as we learn to construct projective resolutions of finite type for certain modules over group rings.
\end{abstract}

\maketitle

\SecRef{Introduction}{Intro}
There are several functorial ways to assign an exact category to a given ring $A$.
Once that is done, Quillen's $K$-theory of the exact category gives an invariant of the ring.

The most common assignment is of the category of finitely generated projective $A$-modules.
Here the admissible monomorphisms are the inclusions of direct summands, and the admissible epimorphisms are the projections onto direct summands.
The class of exact sequences is the class of all short exact sequences, which are split in this case.
The result is called the algebraic $K$-theory of $A$.  

Another useful natural construction can be made when the ring $A$ is Noetherian.
The category is of all finitely generated $A$-modules.
In this case, the admissible monomorphisms are all $A$-linear injections, and the admissible epimorphisms are all surjections.
The class of exact sequences is then the class of all short exact sequences.
The result is called the $G$-theory of $A$. 
While this and related constructions have their applications, the general viewpoint is that $G$-theory is a computable approximation to $K$-theory where most naturally occuring invariants lie.  
There are natural Cartan maps from $K$-theory to $G$-theory given by the relaxation of exactness conditions.
These maps allow the comparison of the two theories.
For example, it is known that when the ring $A$ is regular, the Cartan map is an isomorphism.

In this paper we are interested in the case of $A = R[\Gamma]$, a group ring.
The only class of groups $\Gamma$ for which the group rings are known to be Noetherian, assuming that $R$ is Noetherian, is the class of polycyclic-by-finite groups.
In fact, there is
a conjecture of Philip Hall that only polycyclic-by-finite groups have Noetherian integral group rings.

When $A$ is not Noetherian there is a compromise.  
One can consider the category of all finitely generated $A$-modules as before.
There is a choice of exact sequences which simply consist of short exact sequences where all three modules are finitely generated. 
Of course, this description implies that the admissible monomorphisms are again all injections, but the admissible epimorphisms are only those surjections that have finitely generated kernel.
There is still a Cartan map from $K$-theory to this modified $G$-theory.
Unfortunately, the computational methods are not as good for this $G$-theory.

In this paper we describe a new option.  
It will be developed only for rings $A$ that are group rings $R[\Gamma]$, where $R$ is a Noetherian ring and $\Gamma$ is an arbitrary finitely generated group.
This is the situation of interest for applications in geometric topology.
We will define a new exact structure on a full subcategory $\B (R[\Gamma])$ of finitely generated $R[\Gamma]$-modules.
It is intermediate between the split structure and the modified variant above.
Since it appears to us closer to $K$-theory than to $G$-theory, we will use notation $I (R[\Gamma])$ for the new theory.
To those who believe a certain story about the origin of the traditional notation, we note that ``I'' can also be found between the two extreme letters in the name Grothendieck and is intentionally closer to~``K''.

To make the analogy with $G$-theory complete, in a separate paper \cite{gCbG:13} we show that the new $I$-theory can indeed be effectively computed as the homology of 
$\Gamma$ with coefficients in the usual $G$-theory of $R$ subject to the assumption that $\Gamma$ has a finite classifying space $B\Gamma$.
This is done by showing that an appropriate assembly map is an equivalence.

In the second half of this paper we compute the Cartan map from the $K$-theory of $R[\Gamma]$ to the new $I$-theory.
Using our previous work \cite{gCbG:03}, we prove that it is an equivalence under specific but very general algebraic and geometric assumptions.

\begin{MainThm}
Suppose $R$ is a regular Noetherian ring of finite global homological dimension.
Suppose $\Gamma$ is a discrete group with finite asymptotic dimension and a finite CW-model for $B\Gamma$.
Then the Cartan map $K_n (R[\Gamma]) \to I_n (R[\Gamma])$ is an isomorphism in all dimensions $n \in \mathbb{Z}$.
\end{MainThm}

Our technique of using the coarse geometry of the group $\Gamma$ to construct finite projective resolutions of individual modules from $\B (R[\Gamma])$ should be of independent interest.

\SecRef{Bounded $G$-theory of a Group}{BGOAG}

The $I$-theory is best viewed as the fixed point spectrum of the bounded $G$-theory we developed in \cite{gCbG:00,gCbG:14}.  
The input of the general definition is a proper metric space $X$ and a Noetherian ring $R$.
After a brief review we will specialize to the case of interest in this paper, the group
$\Gamma$ viewed as a word-length metric space with respect to a fixed finite set of generators.

For a proper metric space $X$ and a general ring $R$,
Pedersen and Weibel 
define the category $\mathcal{B} (X,R)$ of \textit{geometric $R$-modules over $X$} in \cite{ePcW:85}.
The objects of this category are locally finite functions $F \colon X \to \Free_{\textit{fg}} (R)$, from points of $X$ to the category of finitely generated free $R$-modules.  We will denote by $F_x$ the module assigned to the point $x$ and the object itself by writing down the collection $\{ F_x \}$.
The \textit{local finiteness} condition requires that for every bounded subset $S \subset X$ the restriction of $f$ to $S$ has finitely many nonzero modules as values.

Let $d$ be the distance function in $X$.  The morphisms $\phi \colon F \to F'$ in $\mathcal{B} (X,R)$ are collections of $R$-linear homomorphisms
$\phi_{x,x'} \colon F_x \longrightarrow F'_{x'}$,
for all $x$ and $x'$ in $X$, with the property
that $\phi_{x,x'}$ is the zero homomorphism whenever $d(x,x') > D$
for some fixed real number $D = D (\phi) \ge 0$.
In this case, we say that $\phi$ is \textit{bounded by} $D$.
The composition of two morphisms $\phi \colon F \to F'$ and $\psi \colon F' \to F''$
is given by the formula
\[
(\psi \circ \phi)_{x,x'} = \sum_{z \in X} \psi_{z,x'} \circ \phi_{x,z}.
\]
The sum in the formula is finite because of the local finiteness of $\{ F'_x \}$.

The category $\mathcal{B} (X,R)$ is an additive category.
The associated $K$-theory is called the \textit{bounded $K$-theory}. 
It is clear that a free action on $X$ by isometries gives an induced free action on the category $\mathcal{B} (X,R)$.
What we need to do next is modify it into an equivariant theory with useful fixed objects.
For details we refer to chapter VI of \cite{gC:95}.

Let $\EGamma$ be the category with the object set $\Gamma$ and the
unique morphism $\mu \colon \gamma_1 \to \gamma_2$ for any pair
$\gamma_1$, $\gamma_2 \in \Gamma$. There is a left $\Gamma$-action
on $\EGamma$ induced by the left multiplication in $\Gamma$.
If $\mathcal{C}$ is a small category with a left $\Gamma$-action,
then the category of functors
$\mathcal{C}_{\Gamma}=\Fun(\EGamma,\mathcal{C})$ is another
category with the $\Gamma$-action given on objects by
the formulas $\gamma(F)(\gamma')=\gamma F (\gamma^{-1} \gamma')$
and $\gamma(F)(\mu)=\gamma F (\gamma^{-1} \mu)$.
It is always
nonequivariantly equivalent to $\mathcal{C}$. The fixed
subcategory $\Fun(\EGamma,\mathcal{C})^{\Gamma} \subset
\mathcal{C}_{\Gamma}$ consists of equivariant functors and
equivariant natural transformations.

Explicitly, when $\mathcal{C} = \mathcal{B} (X,R)$ with the $\Gamma$-action described
above, the objects of $\mathcal{B}_{\Gamma} (X,R)^{\Gamma}$ are the pairs
$(F,\psi)$ where $F \in \mathcal{B} (X,R)$ and $\psi$ is a function on
$\Gamma$ with $\psi (\gamma) \in \Hom (F,\gamma F)$ such that
$\psi(1) = 1$ and $\psi (\gamma_1 \gamma_2) =
\gamma_1 \psi(\gamma_2)  \psi (\gamma_1)$.
These conditions imply that $\psi (\gamma)$ is always an
isomorphism. The set of morphisms $(F,\psi) \to
(F',\psi')$ consists of the morphisms $\phi \colon F \to F'$ in
$\mathcal{B} (X,R)$ such that the squares
\[
\begin{CD}
F @>{\psi (\gamma)}>> \gamma F      
\\
@V{\phi}VV @VV{\gamma \phi}V 
\\
F' @>{\psi' (\gamma)}>> \gamma F'
\end{CD}
\notag
\]
commute for all $\gamma \in \Gamma$. 

A more refined
theory is obtained by replacing $\mathcal{B}_{\Gamma} (X,R)$ with the full
subcategory $\mathcal{B}_{\Gamma, 0} (X,R)$ of functors sending all
morphisms of $\EGamma$ to homomorphisms such that the maps and their inverses are bounded by $0$. 
So $\mathcal{B}_{\Gamma, 0} (X,R)^{\Gamma}$ consists of $(F, \psi)$ with 
$\psi (\gamma)$ of filtration $0$ for all $\gamma \in \Gamma$.

The fixed point category $\mathcal{B}_{\Gamma, 0} (X,R)^{\Gamma}$ is clearly additive.
Now we have designed a spectrum 
$K (X,R)$ as the (nonconnective delooping of) the $K$-theory spectrum of
$\mathcal{B}_{\Gamma, 0} (X,R)$ 
with the desired property
$$K (X,R)^{\Gamma} 
= K \left( \mathcal{B}_{\Gamma, 0} (X,R)^{\Gamma} \right) \simeq  K(R[\Gamma])$$
whenever $\Gamma$ is \textit{geometrically finite} in the sense that there is a finite complex of homotopy type $B\Gamma$. 
We will assume that $\Gamma$ is geometrically finite in the rest of this paper.

Now let $R$ be a Noetherian ring.
At the basic level bounded $G$-theory is locally modeled on the category of finitely generated $R$-modules where the exact sequences are not necessarily split.

We will use the notation $\mathcal{P}(X)$ for the power set of $X$ partially ordered by inclusion. Let $\Mod (R)$ denote the category of left $R$-modules.
If $F$ is a left $R$-module, let $\mathcal{I}(F)$ denote the family of all $R$-submodules of $F$ partially ordered by inclusion.
An $X$\textit{-filtered} $R$-module is a module $F$ together with a functor
$\mathcal{P}(X) \to \mathcal{I}(F)$
such that
the value on $X$ is $F$.
It is \textit{reduced} if $F(\emptyset)=0$.

We will use the notation where
$S[b]$ stands for the
metric $b$-enlargement of $S$ in $X$.
So, in particular, $x[b]$ is the metric ball of radius $b$ centered at $x$.     

An $R$-homomorphism $f \colon F \to F'$ of $X$-filtered modules is \textit{boundedly controlled} if there is a fixed number $b \ge 0$ such that the
image $f (F (S))$ is a submodule of $F' (S [b])$
for all subsets $S$ of $X$.
It is called \textit{boundedly
bicontrolled} if, for some fixed $b \ge 0$, in addition to inclusions of submodules
$f (F (S)) \subset F' (S [b])$,
there are inclusions
$f (F)
\cap F' (S) \subset  f F (S[b])$
for all subsets $S \subset X$.

There is a notion of an admissible exact sequence of filtered modules
which is an exact sequence of $R$-modules with the monic required to be a boundedly bicontrolled monomorphism and the epi required to be a boundedly bicontrolled epimorphism.
For details we refer to section 3 of \cite{gCbG:14}.

\begin{DefRef}{HYUT}
Let $F$ be an $X$-filtered $R$-module.
\begin{enumerate}
\item $F$ is called \textit{lean} or $D$-\textit{lean} if there is a number $D \ge 0$ such that
\[
F(S) \subset \sum_{x \in S} F(x[D])
\]
for every subset $S$ of $X$.
\item $F$ is called \textit{insular} or $d$-\textit{insular} if there is a
number $d \ge 0$ such that
\[
F(S) \cap F(U) \subset F(S[d] \cap U[d])
\]
for every pair of subsets $S$, $U$ of $X$.
\end{enumerate}
\end{DefRef}

The category $\LI (X,R)$ of all lean insular reduced objects and boundedly controlled morphisms is exact with respect to the family of admissible exact sequences.
This category has cokernels but not necessarily kernels.

An $X$-filtered $R$-module $F$ is \textit{locally finitely generated} if $F (S)$ is a finitely generated $R$-module for
every bounded subset $S \subset X$.
The category $\B (X,R)$ is the full
subcategory of $\LI (X,R)$ on the locally finitely generated objects.
It is closed under extensions, so it is an exact category.
There is an exact inclusion
$\mathcal{B}(X,R) \to \B (X,R)$.
The construction of $\mathcal{B}_{\Gamma, 0} (X,R)$ can be replicated here with the result 
$\B_{\Gamma, 0} (X,R)$.
Explicitly, the objects of $\B_{\Gamma, 0} (X,R)$ are functors $\EGamma \to \B (X,R)$ such that all morphisms $\gamma$ map to module isomorphisms $\psi (\gamma)$ in $\B (X,R)$ with the property that both $\psi (\gamma)$ and its inverse are bounded by $0$.

\begin{DefRef}{Gtheory}
Let the spectrum 
$G (X,R)$ be the (nonconnective delooping of) the $K$-theory spectrum of
$\B_{\Gamma, 0} (X,R)$.

There is a controlled equivariant Cartan map $K (X,R) \to G (X,R)$ induced by the exact inclusion, and so is the induced Cartan map
$h \colon K (X,R)^{\Gamma} \to G (X,R)^{\Gamma}$.
\end{DefRef}

We will now specialize to the proper metric on the group $\Gamma$ given as the word metric 
with respect to a finite generating set.
We will then define a certain exact structure on a subcategory 
$\B (R[\Gamma])$ of all finitely generated modules $\Mod_{\textit{fg}} (R[\Gamma])$
and relate it to the exact category $\B_{\Gamma,0} (\Gamma,R)^{\Gamma}$ and to the additive category $\mathcal{B}_{\Gamma,0} (\Gamma,R)^{\Gamma}$.

\begin{DefRef}{AndBack}
Given a finitely generated $R[\Gamma]$-module $F$, fix a finite
generating set $\Sigma$ and define a $\Gamma$-filtration of the
$R$-module $F$ by $F(S) = \langle S \Sigma \rangle_R$, the
$R$-submodule of $F$ generated by $S \Sigma$.  Let $s(F, \Sigma)$
stand for the resulting $\Gamma$-filtered $R$-module.

Conversely, if $F$ is a  $\Gamma$-filtered $R$-module which is
$\Gamma$-\textit{equivariant}, in the sense that $f(\gamma S) = \gamma f(S)$
for all $\gamma \in \Gamma$ and $S \subset \Gamma$, then $F$ has the obvious
$R[\Gamma]$-module structure.
\end{DefRef}

\begin{LemRef}{AlmPull}
Every $R[\Gamma]$-homomorphism $\phi \colon F \to F'$ between finitely generated $R[\Gamma]$-modules
 is bounded as an
$R$-linear homomorphism between the filtered $R$-modules $s (F, \Sigma_F)
\to s (F', \Sigma_{F'})$ with respect to any choice of the finite
generating sets $\Sigma_F$ and $\Sigma_{F'}$.
\end{LemRef}

\begin{proof}
Recall that the word-length metric on
$\Gamma$ is the path metric
induced from the condition that $d (\gamma, \omega \gamma) =
1$ whenever $\gamma \in \Gamma$ and $\omega \in \Omega$.
Consider $x \in F(S) = \langle S \Sigma_F \rangle_R$, then
\[
x = \sum_{s, \sigma} r_{s, \sigma} s \sigma
\]
for a finite collection of pairs $s \in S$, $\sigma \in \Sigma_F$.
Since $F (\{ e \}) = \langle \Sigma_F \rangle_R$ for the identity
element $e$ in $\Gamma$, there is a number $d \ge 0$ such that
$\phi F (\{ e \}) \subset {F'}(e[d])$. Therefore,
\[
\phi (x) = \sum_{s, \sigma} r_{s, \sigma} \phi (s \sigma) =
\sum_{s, \sigma} r_{s, \sigma} s \phi (\sigma) \subset \sum_{s \in
S} s {F'}(e[d]) \subset {F'} (S[d])
\]
because the left translation action by any element $s \in S$ on $e[d]$
in $\Gamma$ is an isometry onto $s[d]$.
\end{proof}

\begin{CorRef}{ChIso}
Given a finitely generated $R[\Gamma]$-module $F$ and two choices
of finite generating sets $\Sigma_1$ and $\Sigma_2$, the filtered
$R$-modules $s (F, \Sigma_1)$ and $s (F, \Sigma_2)$ are isomorphic
as $\Gamma$-filtered $R$-modules.
\end{CorRef}

\begin{proof}
The identity map and its inverse are boundedly controlled as maps
between $s (F, \Sigma_1)$ and $s (F, \Sigma_2)$ by Lemma
\refL{AlmPull}.
\end{proof}

\begin{CorRef}{FilProps}
Finitely generated $R[\Gamma]$-modules $F$ with filtrations $s (F,
\Sigma)$, with respect to arbitrary finite generating sets
$\Sigma$, are locally finite and lean.  If $s (F,
\Sigma)$ is insular and $\Sigma'$ is another finite generating set
then $s (F, \Sigma')$ is also insular.
\end{CorRef}

\begin{proof}
For a finite subset $S$, the submodule $F(S)$ is generated by the
finite set $S\Sigma$. Since $F(x) = \langle x \Sigma \rangle_R$,
\[
F(S) = \sum_{x \in \Sigma} \langle x \Sigma \rangle_R = \langle S
\Sigma \rangle_R,
\]
so $s(F,\Sigma)$ is $0$-lean.  The second
claim follows from Corollary \refC{ChIso}.
\end{proof}

\begin{DefRef}{BR[G]}
We define $\B (R[\Gamma])$ as the full subcategory of $\Mod_{\textit{fg}}
(R[\Gamma])$ on $R$-modules $F$ which are
lean and insular as filtered modules $s (F, \Sigma)$ with respect
to some choice of the finite generating set $\Sigma$.
The objects of $\B (R[\Gamma])$ are referred to as \textit{properly generated} $R[\Gamma]$-modules.

We define $\Btimes (R[\Gamma])$ as the category of
pairs $(F, \Sigma)$ with $F$ in $\B (R[\Gamma])$ and
$\Sigma$ a finite generating set for $F$. The morphisms are
the $R[\Gamma]$-homomorphisms between the modules.
\end{DefRef}

Lemma \refL{AlmPull} shows that the map $s \colon \Btimes (R[\Gamma]) \to \B (\Gamma,R)$
described in Definition \refD{AndBack} is a functor. In fact, it
defines a functor
\[
s_{\Gamma} \colon \Btimes (R[\Gamma]) \longrightarrow
\B_{\Gamma,0} (\Gamma,R)^{\Gamma}
\]
by interpreting $s_{\Gamma} (F, \Sigma) = (F, \psi)$ with $F = s (F,
\Sigma)$ and $\psi (\gamma) \colon F \to \gamma F$ induced from $s
\sigma \mapsto \gamma^{-1} s \sigma$. Since $(\gamma F)(S) =
\langle \gamma^{-1} (S) \Sigma \rangle_R$, it follows that $s_{\Gamma} (F, \Sigma)$ is an object of 
$\B_{\Gamma,0} (\Gamma,R)^{\Gamma}$, and $s$ sends all $R[\Gamma]$-homomorphisms
to $\Gamma$-equivariant homomorphisms.

\begin{LemRef}{ForAppr}
Let $F$ be an object of $\B_{\Gamma,0} (\Gamma,R)^{\Gamma}$, and let $\Sigma$ be
a finite generating set for the $R[\Gamma]$-module $F$. Then the
identity homomorphism $\id \colon s_{\Gamma} (F, \Sigma) \to F$ is
bounded with respect to the induced and the original
filtrations of $F$.
\end{LemRef}

\begin{proof}
If $\Sigma$ is contained in $F(e[d])$, where $e$ is the
identity element in $\Gamma$, then $\gamma \Sigma \subset F(\gamma [d])$ for all $\gamma \in \Gamma$, and
\[
s (F, \Sigma) (S) =
\langle S \Sigma \rangle_R \subset F(S[d])
\]
for all subsets $S
\subset \Gamma$.
\end{proof}

Both functors $s$ and $s_{\Gamma}$ are additive with respect to
the obvious additive structure in $\Btimes (R[\Gamma])$, where
$(F, \Sigma_F) \oplus (G, \Sigma_G) = (F \oplus G, \Sigma_F \times
\Sigma_G)$.

\begin{DefRef}{IUO34}
Let the admissible monomorphisms $\phi \colon
(F, \Sigma_F) \to (F', \Sigma_{F'})$ in $\Btimes (R[\Gamma])$ be the
injections $\phi \colon F \to {F'}$ of $R[\Gamma]$-modules $\phi$
such that $s (\phi) \colon s (F, \Sigma_F) \to s ({F'}, \Sigma_{F'})$ is
a boundedly bicontrolled homomorphism of $\Gamma$-filtered
$R$-modules. This is equivalent to requiring that $s(\phi)$ be an
admissible monomorphism in $\B (\Gamma,R)$. Let the
admissible epimorphisms be the morphisms $\phi$ such that
$s(\phi)$ are admissible epimorphisms in $\B (\Gamma,R)$.
\end{DefRef}

\begin{PropRef}{ItIsEx}
The choice of admissible morphisms defines an exact structure on
$\Btimes (R[\Gamma])$ such that both $s$ and $s_{\Gamma}$ are
exact functors.
\end{PropRef}

\begin{proof}
When checking Quillen's axioms in $\Btimes (R[\Gamma])$, all
required universal constructions are performed in $\B (R[\Gamma])$
with the canonical choices of finite generating sets. In
particular, $\Sigma$ in the pushout $B \cup_A C$ is the image of
the product set $\Sigma_B \times \Sigma_C$ in $B \times C$. The
fact that all candidates for admissible
morphisms are boundedly bicontrolled in $\B (\Gamma,R)$ or
$\B_{\Gamma,0} (\Gamma,R)^{\Gamma}$ are easy to check. 
Exactness of $s$ and $s_{\Gamma}$ is
immediate.
\end{proof}

\begin{DefRef}{ghjjk}
We give $\B (R[\Gamma])$ the minimal exact structure that makes
the forgetful functor
\[
p \colon \Btimes (R[\Gamma]) \to \B (R[\Gamma])
\]
sending $(F, \Sigma)$ to $F$ an exact functor.  In other words, an $R[\Gamma]$-homomorphism
$\phi \colon F \to G$ is an admissible monomorphism or
epimorphism if for some choice of finite generating sets,
$\phi \colon (F, \Sigma_F) \to (G, \Sigma_G)$ is respectively an
admissible monomorphism or epimorphism in $\Btimes (R[\Gamma])$.
\end{DefRef}

Corollary \refC{ChIso} shows that if $\phi \colon F \to G$ is
boundedly bicontrolled as a map of filtered $R$-modules $s (F,
\Sigma_F) \to s (G, \Sigma_G)$ then it is boundedly bicontrolled
with respect to any other choice of finite generating sets, so
this structure is well-defined.

\begin{NotRef}{NewExStr}
The new exact category will be referred to as $\B (R[\Gamma])$,
with the corresponding nonconnective $K$-theory spectrum $I (R[\Gamma])$.
We will use the notation $\Itimes (R[\Gamma])$ for the nonconnective $K$-theory of $\B_{\times} (R[\Gamma])$.
\end{NotRef}

Let $(F,\psi)$ be an object of $\B_{\Gamma,0}
(\Gamma,R)^{\Gamma}$. One may think of $\gamma F \in \B
(\Gamma,R)$, $\gamma \in \Gamma$, as the module $F$ with a new
$\Gamma$-filtration. Now the $R$-module structure $\eta \colon R
\to \End F$ induces an $R[\Gamma]$-module structure $\eta (\psi)
\colon R[\Gamma] \to \End F$ given by
\[
\sum_{\gamma} r_{\gamma} \gamma \mapsto \sum_{\gamma} \eta
(r_{\gamma})
 \psi (\gamma)
\]
since the sums are taken over a finite subset of $\Gamma$. It is
easy to see that sending $(F, \psi)$ to $F$ defines a map
$\pi \colon \B_{\Gamma,0} (\Gamma,R)^{\Gamma} \rightarrow \B
(R[\Gamma])$,
so that $p = \pi s_{\Gamma}$.
Notice, however, that in general $\pi$ is not exact because the identity
homomorphism in Lemma \refL{ForAppr} is not necessarily an
isomorphism.

The
exact functors $p$ and $s_{\Gamma}$ induce maps in
$K$-theory
\[
I (R[\Gamma]) \xleftarrow{\ \, p_{\ast} \ } \Itimes (R[\Gamma])
\xrightarrow{\ s_{\Gamma \ast}\ } G (\Gamma,R)^{\Gamma}.
\]
We claim that both of these maps are weak equivalences.

\begin{PropRef}{Debase}
The functor $p$ induces a weak equivalence
$\Itimes (R[\Gamma]) \simeq I (R[\Gamma])$.
\end{PropRef}

\begin{proof}
This follows from the Approximation Theorem. 
We refer to Theorem 4.5 from \cite{gCbG:00} for the form we use.
The
two categories are saturated, and the derived category of complexes in $\Btimes (R[\Gamma])$ has a
cylinder functor satisfying the cylinder axiom which is
constructed as the canonical homotopy pushout with the canonical
product basis, see section 1 of \cite{rTtT:90}. The first
condition of the Approximation Theorem is clear. For the second
condition, let $(F^{\subdot}_1, \Sigma^{\subdot}_1)$ be a complex
in $\Btimes (R[\Gamma])$ and let $g \colon F^{\subdot}_1 \to
F^{\subdot}_2$ be a chain map in the derived category of $\B (R[\Gamma])$. For each
$R[\Gamma]$-module $F^i_2$ choose any finite generating set
$\Sigma^i_2$, then using $f=g$ and $g'=\id$, we have $g = g'
p(f)$.
\end{proof}

\begin{PropRef}{InjGR}
$s_{\Gamma}$ induces a weak equivalence
$\Itimes (R[\Gamma]) \simeq G (\Gamma,R)^{\Gamma}$.
\end{PropRef}

\begin{proof}
The target category is again saturated and has a cylinder functor
satisfying the cylinder axiom. To check condition (2) of the
Approximation Theorem, let
\[
E^{\subdot} \colon \ 0 \longrightarrow (E^1, \Sigma_1)
\longrightarrow (E^2, \Sigma_2) \longrightarrow\ \dots\
\longrightarrow (E^n, \Sigma_n) \longrightarrow 0
\]
be a complex in $\Btimes (R[\Gamma])$,
\[
(F^{\subdot},\psi.) \colon \ 0 \longrightarrow (F^1, \psi_1)
\xrightarrow{\ f_1\ } (F^2, \psi_2) \xrightarrow{\ f_2\ }\ \dots\
\xrightarrow{\ f_{n-1}\ } (F^n, \psi_n) \longrightarrow 0
\]
be a complex in $\B_{\Gamma,0} (\Gamma,R)^{\Gamma}$, and $g \colon
s'_{\Gamma} (E^{\subdot}) \to (F^{\subdot}, \psi.)$ be a chain
map. Each $F^i$ can be thought of as an $R[\Gamma]$-module, and
there is a chain complex
\[
F^{\subdot} \colon \ 0 \longrightarrow F^1 \xrightarrow{\ f_1\
} F^2 \xrightarrow{\ f_2\ }\ \dots\ \xrightarrow{\ f_{n-1}\ } F^n
\longrightarrow 0
\]
in $\Mod_{\textit{fg}} (R[\Gamma])$. Choose arbitrary finite generating sets
$\Omega_i$ in $F^i$ for all $1 \le i \le n$. Now
\[
\pi_{\Omega} F^{\subdot} \colon \ 0 \longrightarrow (F^1,
\Omega_1) \xrightarrow{\ f_1\ } (F^2, \Omega_2) \xrightarrow{\
f_2\ }\ \dots\ \xrightarrow{\ f_{n-1}\ } (F^n, \Omega_n)
\longrightarrow 0
\]
is a chain complex in $\Btimes (R[\Gamma])$. The chain map $g$
is degree-wise an $R[\Gamma]$-homo\-mor\-phism, so there is a
corresponding chain map $f \colon E^{\subdot} \to \pi_{\Omega}
F^{\subdot}$ which coincides with $g$ on modules. On the other
hand, the degree-wise identity gives a chain map $g' \colon
s'_{\Gamma} (\pi_{\Omega} F^{\subdot}) \to F^{\subdot}$ in
$\B_{\Gamma,0} (\Gamma, R)^{\Gamma}$ by Lemma \refL{ForAppr}. This
$g'$ is a quasi-isomorphism, as required.
\end{proof}

\begin{CorRef}{InjGRreal}
Let $\Gamma$ be a geometrically finite group and $R$ be
a Noetherian ring. There is a weak equivalence
$I (R[\Gamma]) \simeq G (\Gamma ,R)^{\Gamma}$.
\end{CorRef}

\SecRef{Comparison to $K$-theory}{KTC}

We will need to adapt and interpret some of the results from \cite{gCbG:03}.
Those results are stated for equivariant $\Gamma$-modules $F$
over a Noetherian ring $R$ which are locally finite, lean, and satisfy a certain weakening of the insularity property.
Here we restate and prove those results in the category $\B (R[\Gamma])$.

Recall that the objects of $\B (R[\Gamma])$, the \textit{properly generated} $R[\Gamma]$-modules, are the finitely generated $R[\Gamma]$-modules which are
lean and insular as filtered modules $s (F, \Sigma)$ with respect
a choice of the finite generating set $\Sigma$.

\begin{DefRef}{FPMADM}
The group ring $R[\Gamma]$ is 
\textit{weakly regular Noetherian} if every properly generated
$R[\Gamma]$-module has finite
projective dimension.
\end{DefRef}

Our main result in this paper confirms that the class of weakly regular Noetherian rings is surprisingly large and that, as a consequence, the spectrum $I(R[\Gamma])$ is weakly equivalent to the $K$-theory of $R[\Gamma]$ for this class.

\medskip

First we describe the geometric condition on the group $\Gamma$ needed for this result.
The asymptotic dimension is a coarse characteristic of a metric space introduced by M.~Gromov \cite{mG:93}. 

\begin{DefRef}{FAD} 
One considers uniformly bounded families of subsets of a metric space $X$.  Such family $\mathcal{U}$ is called $R$-\textit{disjoint} if for any pair of members $U$ and $U'$ we have $U[R] \cap U' = \emptyset$.  In these terms, $X$ has \textit{asymptotic dimension
at most} $n$ if for any $R > 0$ there exist $n+1$ uniformly bounded $R$-disjoint families which together cover $X$.
It is known that this property is a quasi-isometry invariant.
In particular, when a finitely generated group is given a word metric, having \textit{finite asymptotic dimension} (FAD) as a metric space is an invariant of the group itself.
\end{DefRef}

The class of groups with FAD is very large.
The following families of groups have been verified to have FAD:
Gromov hyperbolic groups \cite{mG:93,jR:05},
one-relator groups \cite{gBaD:05},
virtually polycyclic groups \cite{gBaD:06},
solvable groups with finite rational Hirsch length \cite{aDjS:06},
Coxeter groups \cite{aDtJ:99},
cocompact lattices in connected Lie groups \cite{gCbG:03},
arithmetic groups \cite{lJ:04},
$S$-arithmetic groups \cite{lJ:07,lJ:08},
finitely generated linear group over a field of positive characteristic \cite{eGrTgY:10},
relatively hyperbolic groups with the parabolic subgroup of finite asymptotic dimension \cite{dO:05},
proper isometry groups of finite dimensional $\CAT (0)$ cube complexes, for example B(4)--T(4) small cancellation groups \cite{nW:10},
mapping class groups \cite{mBkBkF:10},
various generalized products of groups from these classes such as fundamental groups of developable complexes of asymptotically finite dimensional groups \cite{gB:02}.

We should mention here that there are known examples of asymptotically infinite dimensional groups, including Thompson's group $F$, Grigorchuk's group, and Gromov's group containing an expander.  The reader can find a discussion of these phenomena in \cite{aDmS:10}.

Our proof of  the main Theorem \refT{Main1} is based on the following
characterization of metric spaces of finite asymptotic dimension and a sequence of lemmas.

\begin{DefRef}{HGT}
A map between metric spaces $\phi \colon (M_1, d_1) \to (M_2,
d_2)$ is a \textit{uniform embedding} if there are two real functions $f$
and $g$ with $\lim_{x \to \infty} f(x) = \infty$ and $\lim_{x \to
\infty} g(x) = \infty$ such that
\[
f(d_1 (x,y)) \le d_2 (\phi(x), \phi(y)) \le g (d_1(x,y))
\]
for all pairs of points $x$, $y$ in $M_1$.
\end{DefRef}

\begin{ThmRefName}{GC}{Dranishnikov \cite{aD:99,aD:00}}
A group $\Gamma$ has finite asymptotic dimension if and only if
there is a uniform embedding of $\Gamma$ in a finite product of
locally finite simplicial trees.
\end{ThmRefName}

\begin{LemRef}{BB}
Let $\Gamma$ be an arbitrary finitely generated group.
Every surjective
$R[\Gamma]$-homomorphism between properly generated $R[\Gamma]$-mo\-du\-les is boundedly
bicontrolled.
\end{LemRef}

\begin{proof}
First, we claim that every $R[\Gamma]$-homomorphism $\phi \colon F \to G$ between
properly generated $R[\Gamma]$-modules is boundedly controlled.
Consider
$z \in F(S)$, then $z = \sum r_i z_i$, where $z_i \in F(x_i [D])$
for some $x_i \in S$, and assuming $F$ is $D$-lean.
Since $\phi$ is an $R[\Gamma]$-homomorphism,
there is a number $b \ge 0$ such that $\phi (z)$ is in
$g(x [D+b])$ for all $z \in f(x [D])$ and all $x \in \Gamma$.
Then
\[
\phi (z) = \sum r_i \phi (z_i) \in \sum g(x_i [D+b])
\subset g(S [D+b]).
\]
If $\phi$ is surjective, and $y \in G(S)$, then $y = \sum r_i y_i$ with $y_i \in G(x_i [D])$, $x_i \in G(S)$, assuming $G$ is also $D$-lean.  Each $G(x_i [D])$ is a finitely
generated $R$-module, so there is a constant $C \ge 0$ and $z_i
\in F (x_i [D+C])$ so that $\phi (z_i) = y_i$.  Now $z = \sum
r_i z_i$ is in $F (S [D+C])$.
\end{proof}

\begin{LemRef}{SCSG}
Let $P$ be a finite product of locally finite simplicial trees, with the
product simplicial metric.
\begin{enumerate}
\item The kernel of a surjective boundedly
bicontrolled homomorphism in $\B (P,R)$ is
also in $\B (P,R)$.

\item If $j \colon M_1 \to M_2$ is a uniform embedding between
proper metric spaces then the $M_2$-filtration $F_{\ast} (S) = F(j^{-1} (S))$ induced from an $M_1$-filtration $F$ is properly generated if and only if $F$ is properly generated.

\item If a group $\Gamma$ has FAD then
the kernel of a surjective $R[\Gamma]$-homomorphism of properly generated
$R[\Gamma]$-modules is properly generated.  In particular, it is finitely
generated.
\end{enumerate}
\end{LemRef}

\begin{proof}
(1) The proof of leanness of the kernel is crucial to the whole argument; it proceeds exactly as that of Lemma 2.4 in \cite{gCbG:03}.
To make this precise, we need to explain one difference in assumptions on filtered modules used in \cite{gCbG:03} and in this paper.
The insularity property in this paper is required for all pairs of subsets.
The property used in the argument in \cite{gCbG:03} is precisely the insularity property applied only to coarsely antithetic pairs of subsets of the metric space.
(Two subsets $S$ and $T$ are called \textit{coarsely antithetic} if for each number $d > 0$ there is $d' > 0$ such that $S[d] \cap T[d] \subset (S \cap T)[d']$.)
So this Lemma has stronger assumptions and a stronger conclusion. 
However, the proof of Lemma 2.4 from \cite{gCbG:03} is valid verbatim in our present setting.

(2) Notice that the fact that $d_2 (j (x), j(y)) \le g (d_1(x,y))$ implies
\[
x[D] \subset j^{-1} (j(x) [g(D)])
\]
for all $D \ge 0$.  Suppose $F$ is $D$-lean, then given $y \in F_{\ast} (S) = F(j^{-1} (S))$ and
\[
y \in \sum_{x \in j^{-1} (S)} F(x[D]),
\]
we have
\begin{align}
y &\in \sum_{x \in j^{-1} (S)} F(j^{-1} (j(x) [g(D)]) \notag \\
&= \sum_{x \in j^{-1} (S)} F_{\ast} (j(x) [g(D)]) \notag \\
&\subset
\sum_{z \in S} F_{\ast} (z [g(D)]). \notag
\end{align}
For insularity, notice that the fact that $f(d_1 (x,y)) \le d_2 (j(x), j(y))$ implies
\[
(j^{-1} (S))[d] \subset j^{-1} (S[f(d)])
\]
for all $d \ge 0$.
If \[
y \in F_{\ast} (S) \cap F_{\ast} (T) = F(j^{-1}(S)) \cap F(j^{-1}(T))
\]
then, assuming $F$ is $d$-insular,
\begin{align}
y &\in F ((j^{-1} (S))[d] \cap (j^{-1} (T))[d]) \notag \\
&\subset F (j^{-1}( S [f(d)] ) \cap j^{-1}( T [f(d)] )) \notag \\
&= F ( j^{-1} ( S [f(d)]  \cap  T [f(d)] )) \notag \\
&= F_{\ast} ( S [f(d)]  \cap  T [f(d)] ). \notag
\end{align}
So $F_{\ast}$ is $g(D)$-lean and $f(d)$-insular.
The sufficiency half of the argument is left to the reader.

(3) Suppose $i \colon \Gamma \to P$ is a uniform embedding given by Dranishnikov's theorem.
Let $\phi \colon F_1 \to F_2$ be the given homomorphism between two
properly generated $R[\Gamma]$-modules.
Now $\phi$ can be thought of as an $R[\Gamma]$-homomorphism between $P$-filtered $R$-modules $F_{1\ast}$ and $F_{2\ast}$ defined by $F_{\ast} (S) = F(i^{-1} (S))$.
By part (2), $F_{\ast}$ is properly generated if and only if $F$ is properly generated.
When $\phi$ is surjective, it is boundedly bicontrolled by Lemma \refL{BB}.  The rest follows from part (1) and the easy fact that lean locally finitely generated $R[\Gamma]$-modules are finitely generated.
\end{proof}

\begin{CorRef}{Conseq}
If the group $\Gamma$ has FAD, the kernel of a surjective $R[\Gamma]$-linear homomorphism between properly generated $R[\Gamma]$-modules is a properly generated $R[\Gamma]$-module.  It is, therefore, finitely generated.
\end{CorRef}

A ring $R$ is said to have \textit{global homological dimension} $\le d$ if every finitely generated left $R$-module $M$ has a resolution
\[
0 \longrightarrow P_{d}
\longrightarrow \ldots \longrightarrow P_{1} \longrightarrow P_0 \longrightarrow M
\longrightarrow 0
\]
where all $P_i$ are finitely generated projective $R$-modules.
If such number $d$ exists then $R$ has \textit{finite global homological
dimension}.

\begin{ThmRef}{Main1}
If $R$ is a Noetherian ring and $\Gamma$ is a finitely generated group with FAD then 
\begin{enumerate}
\item[(1)] all properly generated $R[\Gamma]$-modules have resolutions by finitely generated projective $R[\Gamma]$-modules.
\end{enumerate}
Suppose $R$ is a regular Noetherian ring of finite global homological
dimension and 
suppose $\Gamma$ is a geometrically finite group.
Then we have the following consequences.
\begin{enumerate}
\item[(2)] The group ring $R[\Gamma]$ is weakly regular Noetherian.
\item[(3)] The exact inclusion of bounded categories $\mathcal{C} (\Gamma, R) \to \B (\Gamma, R)$ induces a weak equivalence $K (\Gamma ,R)^{\Gamma} \simeq G (\Gamma ,R)^{\Gamma}$.
\item[(4)] The Cartan map $h \colon K (R[\Gamma]) \to G (R[\Gamma])$ is an equivalence.
\end{enumerate}
\end{ThmRef}

\begin{proof}
(1) Given a properly generated $R[\Gamma]$-module $F$, let $F_0$ be the free
$R[\Gamma]$-module on the finite generating set $\Sigma$ of $F$.
We can view it as a properly generated $\Gamma$-filtered $R$-module with the canonical filtration induced
from the product generating set $\Gamma \times \Sigma$. Then the
surjection of $R$-modules $\pi \colon F_0 \to F$ induced by the map of generating sets $\Gamma \times \Sigma \to \Gamma \Sigma$ given by $(\gamma, \sigma) \mapsto \gamma \sigma$ is boundedly bicontrolled by Lemma \refL{BB}. The kernel $K_1 = \ker (\pi)$ is properly generated by part (3) of Lemma \refL{SCSG}.
Construct
a free finitely generated $R[\Gamma]$-module $F_1$ with a boundedly bicontrolled
projection $\pi_1 \colon F_1 \to K_1$ just as above.
This shows that $F$ is finitely presented as
the quotient of the composition $d_1 = i_1 \pi_1$.
It is an easy fact that the composition of
a bicontrolled epimorphism with any bicontrolled homomorphism is bicontrolled.
Thus $d_1$ is bicontrolled. This construction also inductively gives a resolution
\[
\ldots \longrightarrow F_i
\longrightarrow \ldots \longrightarrow F_1 \longrightarrow F_0 \longrightarrow F
\longrightarrow 0
\]
by finitely
generated free $R[\Gamma]$-modules.

(2) We examine the resolution of $F$ from part (1).
Since $\Gamma$ has a finite $K(\Gamma,1)$ complex, $\Gamma$ belongs to Kropholler's hierarchy $\mathrm{LH} \mathcal{F}$.
By Theorem A of \cite{pK:99},
there is $d \ge 0$ such that the kernel $K_d$ of $F_d \to F_{d-1}$ is isomorphic to a direct summand of a polyelementary module.
When $\Gamma$ is torsion-free, the elementary modules of the form $U \otimes_R R[\Gamma]$, where $U$ is projective over $R$, are themselves projective over $R[\Gamma]$.
So the polyelementary modules are also projective.
Now we have a finite resolution
\[
0 \longrightarrow K_d \longrightarrow F_d
\longrightarrow \ldots \longrightarrow F_1 \longrightarrow F_0 \longrightarrow F
\longrightarrow 0
\]
where $K_d$ is a finitely generated projective $R[\Gamma]$-module.

(3) From Quillen's Resolution theorem \cite{dQ:73} and part (2), applied to the exact inclusion $\mathcal{B}_{\Gamma,0}  (\Gamma ,R)^{\Gamma} \to \B_{\Gamma,0} (\Gamma ,R)^{\Gamma}$, we have $K (\Gamma ,R)^{\Gamma} \simeq G (\Gamma ,R)^{\Gamma}$.

(4) The equivalence $K (\Gamma ,R)^{\Gamma} \simeq K (R[\Gamma])$ was proved for geometrically finite groups in \cite{gC:95}.
The equivalence $K (R[\Gamma]) \simeq \Gnc (R[\Gamma])$ follows from Corollary \refC{InjGRreal}.
\end{proof}

\begin{RemRef}{Lueck}
We would like to point out that our choice of the exact structure in $\B (R[\Gamma])$ is the only one, out of several options pointed out in the introduction, that has the Cartan equivalence in part (4) of Theorem \refT{Main1}.
Recall that, when $R[\Gamma]$ is not Noetherian, the standard compromise for an exact structure on finitely generated $R[\Gamma]$-modules has all injective and surjective $R[\Gamma]$-homomorphisms with
finitely generated cokernels and kernels as admissible morphisms
so that the exact sequences are the usual short exact
sequences.
It can be seen from Remark 2.23 in L\"{u}ck \cite{wL:98} that with the Cartan map to the resulting $G$-theory of $R[\Gamma]$ would not be
injective even in the case when $R$ is the ring of complex number $\mathbb{C}$ and $\Gamma$ is the free group on two generators $F_2$. 
L\"{u}ck shows that the class $[\mathbb{C}]$ is zero in $G_0 (\mathbb{C} [\mathbb{Z}])$.
\end{RemRef}

\SecRef{Examples}{Examples}

One of the reasons we say that $G(R[\Gamma])$ seems to be close to the $K$-theory $K(R[\Gamma])$ is the usual in this area of algebra difficulty in constructing examples of properly generated $R[\Gamma]$-modules that are not projective or even free.
In this section, we want to show some examples in that direction.

The following fact follows from the basic Theorem 2.18 in \cite{gCbG:00}.

\begin{PropRef}{Coker}
Suppose $R$ is a Noetherian ring and $\Gamma$ is a finitely generated group.
The cokernels and images of boundedly bicontrolled homomorphisms between properly generated
$R[\Gamma]$-modules are properly generated.
\end{PropRef}

A boundedly controlled idempotent of a
filtered module is always bicontrolled.  Indeed, if
$\phi \colon F \to F$ is an idempotent so that $\phi^2 = \phi$
then the restriction of $\phi$ to its image is the identity.
Therefore $\phi f(X) \cap f(S) \subset \phi
f(S)$.
So, for example, the cokernels and images of idempotents of finitely generated free $R[\Gamma]$-modules are
examples of properly generated
$R[\Gamma]$-modules.
The results are finitely generated projective $R[\Gamma]$-modules which are often filtered by non-projective $R$-submodules.

The following concrete non-projective examples exploit a similar idea.

\begin{ExRef}{ExInt}
There are systematic constructions of non-free stably free $\mathbb{Z}[\Gamma]$-modules for infinite discrete groups $\Gamma$, cf.~section 17.4 of Johnson \cite{fJ:12}. 
Specifically, there are constructions for geometrically interesting groups by Dunwoody \cite{mD:72}, Berridge--Dunwoody \cite{pBmD:79} for the trefoil group $\langle x,y \mid x^3 = y^2 \rangle$ and similar groups in the work concerned with relation modules. The recent constructions of Harlander--Jensen \cite{jHjJ:06} apply to the Baumslag--Solitar group $G = \langle x,y \mid x y^2 x^{-1} = y^3 \rangle$.
It is also true that $G$ is generated by $x$ and $z=y^4$.
Let $K = R(G,x,z)$ be the kernel of the canonical map $F_2 \to G$ sending the two free generators to $x$ and $z$.
The abelianization $M(G,x,z)=K/[K,K]$ is a $\mathbb{Z}[G]$-module, where the action of $G$ is by conjugation, called the relation module associated with the generating set $\{ x,z \}$. Now Harlander and Jensen show that $M$ is non free but $M \oplus \mathbb{Z}[G] \cong \mathbb{Z}[G]^2$.
The projection onto $M$ is an interesting idempotent $e$ of $\mathbb{Z}[G]^2$.
Let $e_n \colon \mathbb{Z}_n [G]^2 \to \mathbb{Z}_n [G]^2$ be the reduction modulo $n > 1$. The image $M_n$ of $e_n$ is not flat over $\mathbb{Z} [G]$, so it is not projective.
However, $e_n$ is a bicontrolled endomorphism of a properly generated $\mathbb{Z}[G]$-module $\mathbb{Z}_n [G]^2$, so $M_n$ is properly generated.  
\end{ExRef}

\end{document}